\documentclass[11pt]{article}
\usepackage{amscd}
\usepackage{amsfonts}
\usepackage{amsmath}
\usepackage{amssymb}
\usepackage{amsthm}
\usepackage{bbm}
\usepackage{CJK}
\usepackage{fancyhdr}
\usepackage{graphicx}
\usepackage{hyperref}
\usepackage{indentfirst}
\usepackage{latexsym}
\usepackage{mathrsfs}
\usepackage{xypic}

\usepackage[top=1in,bottom=1in,left=1.25in,right=1.25in]{geometry}
\newtheorem{theorem}{Theorem}[section]
\newtheorem{lemma}[theorem]{Lemma}
\newtheorem{definition}[theorem]{Definition}
\newtheorem{proposition}[theorem]{Proposition}
\newtheorem{example}[theorem]{Example}

\newtheorem{corollary}[theorem]{Corollary}
\newtheorem{remark}[theorem]{Remark}

\usepackage[top=1in,bottom=1in,left=1.25in,right=1.25in]{geometry}
\def\<{\langle}
\def\>{\rangle}
\def\a{\alpha}
\def\b{\beta}

\def\g{\gamma}

\def\o{\otimes}

\date{}
\begin{document}
\renewcommand{\baselinestretch}{1.2}
\renewcommand{\arraystretch}{1.0}
\title{\bf The construction  of  Hom left-symmetric conformal bialgebras}
\author{{\bf Shuangjian Guo$^{1}$, Xiaohui Zhang$^{2}$,  Shengxiang Wang$^{3}$\footnote
        { Corresponding author(Shengxiang Wang):~~wangsx-math@163.com} }\\
{\small 1. School of Mathematics and Statistics, Guizhou University of Finance and Economics} \\
{\small  Guiyang  550025, P. R. of China} \\
{\small 2.  School of Mathematical Sciences, Qufu Normal University}\\
{\small Qufu  273165, P. R. of China}\\
{\small 3.~ School of Mathematics and Finance, Chuzhou University}\\
 {\small   Chuzhou 239000,  P. R. of China}}
 \maketitle
\begin{center}
\begin{minipage}{13.cm}

{\bf \begin{center} ABSTRACT \end{center}}
 In this paper,  we first introduce the notion of Hom-left-symmetric conformal bialgebras and show some nontrivial examples.
 Also,  we present construction methods of matched
pairs of Hom-Lie conformal algebras and Hom-left-symmetric conformal algebras. Finally,  we
prove that a finite Hom-left-symmetric conformal bialgebra  is free as a $\mathbb{C}[\partial]$-module is
equivalent to a Hom-parak\"{a}hler Lie conformal algebra. In particular, we investigate the coboundary Hom-left-symmetric conformal  bialgebras.

{\bf Key words}:  Hom-Lie conformal algebra, Hom-left-symmetric conformal coalgebra,  Hom-left-symmetric conformal bialgebra,
Hom-parak\"{a}hler Lie conformal algebra.

 {\bf 2010 Mathematics Subject Classification:} 17A30, 17B45, 17D25, 17B81
 \end{minipage}
 \end{center}
 \normalsize\vskip1cm

\section*{INTRODUCTION}
\def\theequation{0. \arabic{equation}}
\setcounter{equation} {0}

Left-symmetric algebras are a class of Lie-admissible algebras
whose commutators are Lie algebras. They arose from the study of convex
homogeneous cones, affine manifolds and affine structures on Lie groups.
 Burde gave a survey about left-symmetric algebras
which play an important role in many fields in mathematics and mathematical
physics such as vector fields, rooted tree algebras, words in two letters, vertex
algebras, operad theory, deformation complexes of algebras, convex homogeneous
cones, affine manifolds, left-invariant affine structures on Lie groups in \cite{Burde2006}.
 Furthermore, a theory of left-symmetric bialgebras was developed by Bai in \cite{Bai2008},
he proved that a left-symmetric bialgebra is equivalent to a
parak¡§ahler Lie algebra which is the Lie algebra of a Lie group $G$ with a $G$-invariant
parak\"{a}hler structure, studied coboundary left-symmetric bialgebras
and obtained an analog of the classical Yang-Baxter equation and called it $S$-equation.

Lie conformal algebras were introduced by Kac in \cite{Kac98},
he gave an axiomatic description of the singular part of the operator product expansion of chiral fields in conformal field theory.
 It is an useful tool to study vertex algebras and has many applications in the theory of Lie algebras.
 Moreover, Lie conformal algebras  have close connections to Hamiltonian formalism in the theory of nonlinear evolution equations.
 Lie conformal algebras were widely studied in the following aspects:
  the structure theory \cite{D'Andrea1998}, representation theory \cite{Cheng1997}, \cite{Cheng1998} and cohomology
theory \cite{Bakalov1999} of finite Lie conformal algebras.
Later, Liberati
in \cite{L08} introduced a conformal analog of Lie bialgebras including the conformal classical Yang-Baxter equation,
 the conformal Manin triples and conformal Drinfeld's double.
  As a generalization of \cite{Bai2008}, Hong and Li introduced the definition of left-symmetric conformal algebra in \cite{H2015}
  and developed a conformal theory of left-symmetric bialgebras  in \cite{Hong2015}.

  Recently, the Hom-Lie conformal algebra was introduced
and studied in \cite{Yuan14}, where it was proved that a Hom-Lie conformal algebra is equivalent to a Hom-Gel'fand-Dorfman bialgebra.
From then on, similar generalizations
of certain algebraic structures became a very popular subject.
In \cite{Sheng2014}, Sheng and Bai introduced a new definition of Hom-Lie algebras and studied their properties.
 In \cite{Sun2017}, Sun and Li  proved that a Hom-left-symmetric bialgebra is equivalent to a
Hom-parak¡§ahler Lie algebra, extending the result  given in \cite{Bai2008}.
Recently,  Zhao, Yuan and Chen  developed the cohomology theory of Hom-Lie conformal algebras and discuss
some applications to the study of deformations of regular Hom-Lie conformal
algebras. Also, they  introduced the notion of  derivations of multiplicative Hom-Lie conformal
algebras and study their properties in \cite{Zhao2018}.

Motivated by these results, this paper is organized as follows. In Section 2,  we introduce the definition of a parak\"{a}hler   Hom-Lie conformal algebra, and study the parak\"{a}hler   Hom-Lie conformal algebra in terms of Hom-left-symmetric conformal algebras. In Section 3,  we introduce the matched pairs of Hom-Lie conformal algebras and Hom-left-symmetric conformal algebras. Moreover, we study the relations
between them. In Section 4, we introduce the notion of  a Hom-left-symmetric conformal  bialgebra, which is equivalent to the parak\"{a}hler Hom-Lie conformal algebra. In particular, we investigate the coboundary Hom-left-symmetric conformal  bialgebras.

\section{Preliminaries}
\def\theequation{\arabic{section}.\arabic{equation}}
\setcounter{equation} {0}

Throughout the paper, all algebraic systems are supposed to be over a field $\mathbb{C}$,
 and denote by $\mathbb{Z}_{+}$ the set of all nonnegative integers and by $\mathbb{Z}$ the set of all integers.

In this section we recall some basic definitions and results related to our paper from \cite{L08} and \cite{Yuan14}.

\begin{definition}
A Hom-Lie conformal algebra $(R, \a)$ is a   $\mathbb{C}[\partial]$-module equipped with an even linear endomorphism
$\a$ such that $\a\partial=\partial\a$, and  a $\mathbb{C}$-linear map
\begin{eqnarray*}
R\o R\rightarrow \mathbb{C}[\lambda]\o R,  ~~~a\o b\mapsto [a_\lambda b]
\end{eqnarray*}
satisfying the following axioms:
\begin{eqnarray*}
&&[\partial a_{\lambda}b] =-\lambda[a_\lambda b],[ a_{\lambda}\partial b] =(\partial+\lambda)[a_\lambda b],\\
&&[a_{\lambda}b]=-[b_{-\lambda-\partial}a],\\
&&[\a(a)_\lambda[b_\mu c]]=[[a_{\lambda}b]_{\lambda+\mu}\a(c)]+[\a(b)_{\mu}[a_{\lambda}c]],
\end{eqnarray*}
  for any  $a,b,c\in R.$
\end{definition}
A Hom-Lie conformal algebra $(R, \a)$  is called multiplicative if $\a$ is an algebra endomorphism, i.e., $\a([a_\lambda b])=[\a(a)_{\lambda}\a(b)]$ for any $a,b\in R$. In particular, if $\a$ is an algebra isomorphism, then $(R, \a)$  is regular.

\begin{definition}
Let $(U, \a)$ and $(V, \b)$ be two $\mathbb{C}[\partial]$-modules. A  conformal linear map from
$U$ to $V$ is a $\mathbb{C}$-linear map $a:U\rightarrow \mathbb{C}[\lambda]\o V$, denoted by $a_\lambda: U\rightarrow V$,
such that
 $$[\partial,a_{\lambda}]=-\lambda a_{\lambda}, \partial_U\a=\a\partial_U,  \partial_V\b=\b\partial_V, a_{\lambda}\a=\b a_{\lambda}.$$
 Denote the $\mathbb{C}$-vector space of all such maps by $Chom(U,V)$. It has a  canonical structure of a $\mathbb{C}[\partial]$-module:
$
(\partial a)_{\lambda}=-\lambda a_{\lambda}.
$
\end{definition}
Define the conformal dual of a $\mathbb{C}[\partial]$-module $U$ as $U^{\ast c}= Chom(U,C)$, where $\mathbb{C}$ is
viewed as the trivial $\mathbb{C}[\partial]$-module, that is
\begin{eqnarray*}
U^{\ast c}= \{a:U\rightarrow \mathbb{C}[\lambda]|\mbox{$a$ is $\mathbb{C}$-linear and $a_{\lambda}(\partial b)=\lambda a_{\lambda}b$}\}.
\end{eqnarray*}
Set $Cend(V)=Chom(V,V)$ and assume that $(V, \b)$
is a finite $\mathbb{C}[\partial]$-module, then the $\mathbb{C}[\partial]$-module $Cend(V)$ has a canonical structure defined by
\begin{eqnarray*}
 (a_{\lambda}b)_{\mu}\b(v)=\psi(a_{\lambda})(b_{\mu-\lambda}v).
\end{eqnarray*}
 for any $a,b\in Cend(V), v\in V$ and $\psi:Cend(V)\rightarrow Cend(V).$
Therefore, $(gc(V):=Chom(V,V), \psi)$ has a Hom-Lie conformal algebra structure defined by
\begin{eqnarray*}
[a_{\lambda}b]\b(v)=\psi(a_{\lambda})(b_{\mu-\lambda}v)-\psi(b_{\mu-\lambda})(a_{\lambda}v).
\end{eqnarray*}
Here $(gc(V), \psi)$ is called the general Hom-Lie conformal algebra of $V$.
\section{Parak\"{a}hler   Hom-Lie conformal algebras }
\def\theequation{\arabic{section}. \arabic{equation}}
\setcounter{equation} {0}

In this section, we introduce the definition of a parak\"{a}hler   Hom-Lie conformal algebra
 and study the parak\"{a}hler   Hom-Lie conformal algebra in terms of Hom-left-symmetric conformal algebras.

\begin{definition}
A  Hom-left-symmetric conformal algebra   $(A, \a)$  is a $\mathbb{C}[\partial]$-module equipped with an even linear endomorphism
$\a$ such that $\a\partial=\partial\a$ and  a linear map
\begin{eqnarray*}
A\o A\rightarrow \mathbb{C}[\lambda]\o A,  ~~~a\o b\mapsto a_\lambda b
\end{eqnarray*}
satisfying the following axioms
\begin{eqnarray}
&&(\partial a_{\lambda}b) =-\lambda(a_\lambda b), a_{\lambda}\partial b =(\partial+\lambda)(a_\lambda b),\a(a_\lambda b)=\a(a)_{\lambda}\a(b),\\
&&(a_\lambda b)_{\lambda+\mu} \a(c)-\a(a)_{\lambda}(b_\mu c)=(b_\mu a)_{\lambda+\mu} \a(c)-\a(b)_{\mu}(a_\lambda c),
\end{eqnarray}
for any $a,b,c\in A.$
\end{definition}
A Hom-left-symmetric conformal algebra $(A, \a)$ is called finite if $A$ is a finitely generated $\mathbb{C}[\partial]$-module.
The rank of $A$ is its rank as a $\mathbb{C}[\partial]$-module.

\begin{example}
Let $(A, \a)$ be a Hom-left-symmetric algebra. Then we can naturally define a  Hom-left-symmetric conformal algebra $Cur A=\mathbb{C}[\partial]\o A$ with the $\lambda$-product
\begin{eqnarray*}
a_\lambda b=ab, ~~~~\forall a,b\in A.
\end{eqnarray*}
$Cur A$ is called the current Hom-left-symmetric conformal algebra. $Coeff Cur A= A\o \mathbb{C}[t, t^{-1}]$ with the left-symmetric multiplication
\begin{eqnarray*}
at^{m}bt^{n}=(ab)t^{m+n}, \a(f(\partial)a)=f(\partial)\a(a).
\end{eqnarray*}
\end{example}

\begin{definition}
A Hom-Novikov conformal algebra $(A, \a)$ is a Hom-left-symmetric conformal algebra satisfying
\begin{eqnarray}
(a_\lambda b)_{\lambda+\mu}\a(c)=(a_\lambda c)_{-\mu-\partial} \a(b),
\end{eqnarray}
for any $a,b,c\in R.$
\end{definition}

\begin{remark}
By (2.1), the following equalities always hold.
\begin{eqnarray*}
&& (a_{-\lambda-\partial}b)_{\lambda+\mu}\a(c)=(a_{\mu}b)_{\lambda+\mu}\a(c),\\
&& \a(a)_{\mu}(b_{-\lambda-\partial}c)=\a(a)_{\mu}(b_{-\lambda-\partial-\mu}c).
\end{eqnarray*}
Moreover, by the equalities above and (2.3), we have
\begin{eqnarray*}
(a_\lambda b)_{-\mu-\partial} \a(c)-\a(a)_{\lambda}(b_{-\mu-\partial} c)=(b_{-\lambda-\partial} a)_{-\mu-\partial} \a(c)-\a(b)_{-\mu-\partial-\lambda}(a_\lambda c), \forall a,b,c\in R.
\end{eqnarray*}
\end{remark}

The following result is easily shown.

\begin{proposition}
Let $(A,\a)$ be a  Hom-left-symmetric conformal algebra.    Then the $\lambda$-bracket
\begin{eqnarray*}
[a_{\lambda}b] =a_{\lambda}b-b_{-\lambda-\partial}a, \forall a, b\in A,
\end{eqnarray*}
define a Hom-Lie conformal algebra $(\mathfrak{g}(A), \a)$, which is called the sub-adjacent  Hom-Lie conformal algebra of $(A,\a)$.
 $(A,\a)$ is also called the compatible Hom-left-symmetric conformal algebra structure on the Hom-Lie conformal algebra $(\mathfrak{g}(A), \a)$,
\end{proposition}

\begin{example}
Given a Hom-Lie conformal algebra $R=\mathbb{C}[\partial]a\oplus \mathbb{C}[\partial]b$  with a $\lambda$-bracket:
\begin{eqnarray*}
&& \a(a)=\lambda a,  \a(b)=b,\\
&&[a_\lambda a] =0, [a_\lambda b] =d(\lambda)b, [b_\lambda b]=0,
\end{eqnarray*}
where $d(\lambda)\in \mathbb{C}[\lambda]$.
Then we have a compatible Hom-left-symmetric conformal algebra $(\mathfrak{g}(A), \a)$ with a $\lambda$-bracket:
\begin{eqnarray*}
a_{\lambda}a=0, a_{\lambda}b=d(\lambda)b, b_\lambda a=0,b_\lambda b=0.
\end{eqnarray*}
\end{example}

Let $V$ be a $\mathbb{C}[\partial]$-module, recall from \cite{L08}, a conformal bilinear form on $V$ is a $\mathbb{C}$-bilinear map $\omega_{\lambda}: V\o V\rightarrow \mathbb{C}[\lambda]$ such that
\begin{eqnarray*}
\omega(\partial v, w)_{\lambda}=-\lambda\omega(v, w)_{\lambda}=-\omega(v, \partial w)_{\lambda},
\end{eqnarray*}
 for any $v,w\in V$.
The conformal  bilinear form  is called skew-symmetric if $\omega(v, w)_{\lambda}=-\omega(w, v)_{-\lambda}$ for any $v,w\in V$.
\medskip

Suppose that there is a conformal  bilinear form on a $\mathbb{C}[\partial]$-module $V$. Then we have a $\mathbb{C}[\partial]$-module homomorphism $T: V\rightarrow V^{\ast c},  (T_v)_\lambda w=\omega(v,w)_\lambda$.
\medskip

A conformal bilinear form $\omega_{\lambda}$ is called non-degenerate if $T$ is isomorphism of $\mathbb{C}[\partial]$-modules.

\begin{definition}
A  multiplicative  Hom-Lie conformal algebra  $(R, \a)$ is called a symplectic Hom-Lie conformal algebra if there is a nondegenerate skew-symmetric conformal bilinear form $\omega_{\lambda}$  on $R$ such that
\begin{eqnarray*}
\omega([a_\lambda b], \a(c))_{\mu}+\omega([b_{\mu-\partial} c], \a(a))_{-\lambda}+\omega([c_{-\mu} a], \a(b))_{\lambda-\mu}=0,
\end{eqnarray*}
 for any $a,b,c\in R.$
\end{definition}

\begin{example}
Let $R=\mathbb{C}[\partial]L\oplus \mathbb{C}[\partial]E$. Then $(R,\a)$ is a Hom-Lie conformal algebra with a $\lambda$-bracket:
\begin{eqnarray*}
&&  \a(L)=f(\partial)L, \a(E)=g(\partial)E,\\
&&[L_\lambda L] =(\partial+2\lambda)E, [L_\lambda E]=[E_\lambda E]=0.
\end{eqnarray*}
Then $(R,\a)$ is a symplectic Hom-Lie conformal algebra with a conformal bilinear form $\omega_{\lambda}$:
\begin{eqnarray*}
\omega(L,L)_{\lambda}=0, \omega(L,E)_{\lambda}=1, \omega(E,E)_{\lambda}=0.
\end{eqnarray*}
\end{example}

\begin{example}
Let $(g,\a, \omega)$ be a symplectic Hom-Lie algebra. Then $(Cur g, \omega_\lambda)$ is a symplectic Hom-Lie conformal algebra
with a nondegenerate skew-symmetric conformal bilinear form $\omega_\lambda$:
\begin{eqnarray*}
\omega(p(\partial)a, q(\partial)b)_{\lambda} =p(-\lambda)q(\lambda)\omega(a, b), \forall p(\partial), q(\partial) \in \mathbb{C}[\partial], a,b\in g.
\end{eqnarray*}
\end{example}

\begin{proposition}
Let $(R,\a, \omega_{\lambda})$ be a symplectic Hom-Lie conformal algebra. Then there exists a compatible Hom-left-symmetric conformal algebra structure on $R$ with
\begin{eqnarray}
\omega(a_\lambda b, \a(c))_{\mu}=-\omega(\a(b), [a_\lambda c])_{\mu-\lambda}, \forall a,b,c\in R.
\end{eqnarray}
\end{proposition}

{\bf Proof. }  For any $a,b,c\in R$, we have
\begin{eqnarray*}
&& \omega([a_\lambda b], \a(c))_{\mu}\\
&=& -\omega([b_{\mu-\partial} c], \a(a))_{-\lambda}-\omega([c_{-\mu} a], \a(b))_{\lambda-\mu}\\
&=& \omega([a_{\lambda} c], \a(b))_{\lambda-\mu}-\omega([b_{-\lambda+\mu} c], \a(a))_{-\lambda}\\
&=& -\omega(\a(b), [a_{\lambda} c])_{-\lambda+\mu}+\omega(\a(a), [b_{-\lambda+\mu} c])_{\lambda}\\
&=& \omega(a_\lambda b-b_{-\lambda-\partial}, \a(c))_{\mu}.
\end{eqnarray*}
Since $\omega_{\lambda}$ is non-degenerate, we have $[a_\lambda b]=a_\lambda b-b_{-\lambda-\partial}$.
Also, we have
\begin{eqnarray*}
&&\omega((a_\lambda b)_{\lambda+\mu} \a(c)-\a(a)_{\lambda}(b_\mu c)-(b_\mu a)_{\lambda+\mu} \a(c)+\a(b)_{\mu}(a_\lambda c), \a^2(d))\\
&=& -\omega(\a^2(c), [(a_\lambda b)_{\lambda+\mu}\a(d)])_{\nu-\lambda-\mu}+\omega(\a(b_\mu c), [a_\lambda d])_{\nu-\lambda}\\
&& +\omega(\a^2(c),[(b_\mu a)_{\lambda+\mu}\a(d) ])_{\nu-\lambda-\mu}-\omega(\a(a_\mu c), [b_\mu d])_{\nu-\mu}\\
&=& -\omega(\a^2(c), [(a_\lambda b-b_{-\lambda-\partial})_{\lambda+\mu}\a(d)])_{\nu-\lambda-\mu}-\omega(\a(c), [\a(b)_\mu [a_\lambda d]])_{\nu-\lambda-\mu}\\
&& +\omega(\a^2(c), [\a(a)_{\lambda}[b_{\mu}d]])_{\nu-\lambda-\mu}\\
&=& -\omega(\a^2(c), [(a_\lambda b-b_{-\lambda-\partial})_{\lambda+\mu}\a(d)]+[\a(b)_\mu [a_\lambda d]]+[\a(a)_{\lambda}[b_{\mu}d]])=0.
\end{eqnarray*}
 Therefore,  $(R,\a)$ is a Hom-left-symmetric conformal algebra.  \hfill $\square$

\begin{definition}
Let $(R,\a)$ be a  Hom-Lie conformal algebra. Then $(R,R_0,R_1, \omega_\lambda, \a)$ is called a parak\"{a}hler   Hom-Lie conformal algebra,  if $R_0$  and $R_1$  are two  Hom-Lie conformal subalgebras of $R$ such that $R=R_0\oplus R_1$  as a $\mathbb{C}[\partial]$-module,  and $\omega(R_i,R_i)_\lambda=0$ for $i=0,1$.
\end{definition}

\begin{example}
Let $(g,g_0,g_1, \omega, \a)$ be a parak\"{a}hler   Hom-Lie  algebra. Then $(Cur g, Cur g_0,$\\$ Cur g_1, \omega_\lambda, \a)$ is also a parak\"{a}hler   Hom-Lie conformal algebra with the following structures:
\begin{eqnarray*}
\omega(p(\partial)a, q(\partial)b)_{\lambda} =p(-\lambda)q(\lambda)\omega(a, b),
\end{eqnarray*}
 for any $p(\partial), q(\partial) \in \mathbb{C}[\partial]$ and $a,b\in g.$
\end{example}

\begin{proposition}
Let $(R,R_0,R_1, \omega_\lambda, \a)$ be a parak\"{a}hler   Hom-Lie conformal algebra. Then (2.4) defines a compatible Hom-left-symmetric conformal  algebra structure on $(R, \a)$.
\end{proposition}

{\bf Proof.} Directly from Proposition 2.9. \hfill $\square$

\section{ Matched pairs of  Hom-left-symmetric conformal algebras }
\def\theequation{\arabic{section}. \arabic{equation}}
\setcounter{equation} {0}

In this section, we study the matched pairs of Hom-Lie conformal algebras and
Hom-left-symmetric conformal algebras.

\begin{proposition}$^{\cite{Zhao2018}}$
Let $(R,\a)$ be a Hom-Lie conformal algebra and $(M, \b)$  a module over $R$. Then the  $\mathbb{C}[\partial]$-module $(R\oplus M, \a+\b)$ can be endowed with a Hom-Lie conformal algebra structure as follows:
\begin{eqnarray*}
&& [(a+u)_{\lambda}(b+v)]=[a_{\lambda}b]+a_{\lambda}v-b_{-\lambda-\partial} u,\\
&& (\a+\b)(a+u)=\a(a)+\b(u),
\end{eqnarray*}
for any $a,b\in A$  and $u,v\in M$. Denote this Hom-Lie conformal algebra by $R\ltimes M$.

\end{proposition}

\begin{proposition}
Let $(R, \a)$ and $(R', \a')$ be two Hom-Lie conformal algebras. Suppose that $\rho:R\rightarrow gc(R')$  and $\sigma: R'\rightarrow gc(R)$ are two representations satisfying  the following conditions:
\begin{eqnarray}
&& \rho(\a(x)_{\lambda}[a_{\mu}b]-[(\rho(x)_{\lambda}a)_{\lambda+\mu}\a'(b)]-[\a'(a)_{\mu}(\rho(x)_{\lambda}b)]+
\rho(\sigma(a)_{-\lambda-\partial}x)_{\lambda+\mu}\a'(b)\nonumber\\
&& -\rho(\sigma(b)_{-\lambda-\partial}x)_{-\mu-\partial}\a'(a)=0,\\
&& \sigma(\a'(a)_{-\lambda-\mu-\partial}[x_{\lambda}y]-[\a(x)_{\lambda}(\sigma(a)_{-\mu-\partial} y)]+[\a(y)_{\mu}(\sigma(a)_{-\lambda-\partial}x)]\nonumber\\
&&+\sigma(\rho(x)_{\lambda}a)_{-\mu-\partial}\a(y) -\sigma(\rho(y)_{\mu}a)_{-\lambda-\partial}\a(x)=0,
\end{eqnarray}
for any $x,y\in R$ and $a,b\in R'$. Then their is a Hom-Lie conformal algebra structure on $\mathbb{C}[\partial]$-module $(R\oplus R', \a+\a')$:
\begin{eqnarray}
&&[(x+a)_{\lambda}(y+b)]=[x_\lambda y] +\sigma(a)_{\lambda}y-\sigma(b)_{-\lambda-\partial}x\nonumber \\
&& +[a_\lambda b] +\rho(x)_{\lambda}b-\rho(y)_{-\lambda-\partial}a,\\
&& (\a+\a')(x+a)=\a(x)+\a'(a),
\end{eqnarray}
for any $x,y\in R$ and $a,b\in R'$. We denote this Hom-Lie conformal algebra by $R\bowtie R'$.
If $(R,R',\a,a', \rho, \sigma)$ satisfying the above conditions, then we call it a matched pair of Hom-Lie conformal algebras.
\end{proposition}

{\bf Proof.}
For any $x,y, z\in R$ and $a,b,c\in R'$, we get
\begin{eqnarray*}
&&R(x+a,  y+b,z+c)_{(\a+\a')}=[(\a+\a')(x+a)_{\lambda}[(y+b)_{\mu}(z+c)]]\\
&&- [[(x+a)_{\lambda}(y+b)]_{\lambda+\mu}(\a+\a')(z+c)]]-[(\a+\a')(y+b)_{\mu} [(x+a)_{\lambda}(z+c)]].
\end{eqnarray*}
By direct computation, $R(x,y,z)=0$ if and only if $(R,\a)$ is a Hom-Lie conformal algebra, and $R(a,b,c)=0$ if and only if $(R',\a')$ is a Hom-Lie conformal algebra. $R(x,y,a)=0$ or $R(x,a,y)=0$ or $R(a,x,y)=0$ if and only if $\rho$ is a representation of $(R, \a)$ and (3.3) holds,  $R(x,a,b)=0$ or $R(a,x,b)=0$ or $R(a,b,x)=0$ if and only if $\sigma$ is a representation of $(R', \a')$ and (3.4) holds.   And this the proof. \hfill $\square$

\begin{definition}
Let $(A,\a)$ be a  Hom-left-symmetric conformal algebra and $(M,\b)$ a module. Then we call $(M,\b)$  a $\mathbb{C}[\partial]$-module if there is an even linear endomorphism
$\a$ such that $\b\partial=\partial\b$  and two $\mathbb{C}$-bilinear maps
$$A\o M\rightarrow \mathbb{C}[\lambda]\o M, a\o m\mapsto a_{\lambda}m,~ M \o  A \rightarrow \mathbb{C}[\lambda]\o M, m\o a\mapsto m_{\lambda}a$$
 satisfying
\begin{eqnarray*}
&& (\partial a)_{\lambda} m=-\lambda a_{\lambda}m, (\partial m)_{\lambda} a=-\lambda m_{\lambda}a, \\
&& (a_{\lambda}b)_{\lambda+\mu} \b(m)-\a(a)_{\lambda}(b_\mu m)=(b_\mu a)_{\lambda+\mu}\b(m)-\a(b)_{\mu}(a_\lambda m),\\
&& (a_{\lambda}m)_{\lambda+\mu} \a(b)-\a(a)_{\lambda}(m_\mu b)=(m_\mu a)_{\lambda+\mu}\a(b)-\b(m)_{\mu}(a_\lambda b),
\end{eqnarray*}
for any $a,b\in A$ and $m\in M$.
\end{definition}

\begin{definition}
 Let $(A, \a)$ be a Hom-left-symmetric conformal algebra and $(M, \b)$ a  $\mathbb{C}[\partial]$-module.
Assume $l_A,r_A: A\rightarrow Cend(M)$ are two $\mathbb{C}[\partial]$-module homomorphisms. $(M, l_A, r_A,  \b)$ is called a module if
 \begin{eqnarray}
&&\b(l_A(a)m)=l_A(\a(a))\b(m), \b(r_A(a)m)=r_A(\a(a))\b(m),\\
&& l_A(\partial a)_{\lambda}v=-\lambda l_A(a)_{\lambda}v, r_A(a)_{\lambda}(\partial v)=-\lambda r_A(a)_{\lambda}(v),\nonumber\\
&& l_A(a_\lambda b)_{\lambda+\mu}\b(v)-l_A(\a(a))_\lambda(l_A(b)v)=l_A(b_{\mu}a)_{\lambda+\mu}\b(v)-l_A(\a(b))_{\mu}  (l_A(a)_{\lambda}v),~~~\\
&& r_A(\a(b))_{-\lambda-\mu-\partial}(l_A(a)_\lambda v)-l_A(\a(a))_\lambda(r_A(b)_{-\mu-\partial}v)\nonumber\\
&& =r_A(\a(b))_{-\lambda-\partial\mu}(r_A(a)_{\lambda}v)-r_A(a_\lambda b)_{-\mu-\partial} \b(v),
 \end{eqnarray}
  for any $a,b\in A$ and $m\in M$.
 \end{definition}

\begin{proposition}
 Let $(A, \a)$ be a Hom-left-symmetric conformal algebra and $(M, \b)$ a  $\mathbb{C}[\partial]$-module.
 Assume $l_A,r_A: A\rightarrow Cend(M)$ are two $\mathbb{C}[\partial]$-module homomorphisms. Then $(M, l_A, r_A,  \b)$ is  a module if and only if $(A\oplus M, \a+\b)$ is a Hom-left-symmetric conformal algebra with the $\lambda$-product
 \begin{eqnarray*}
&&(a+u)_{\lambda}(b+v)=a_\lambda b+l_A(a)_\lambda v+r(b)_{-\lambda-\partial} u,\\
&&(\a+\b)(a+u)=\a(a)+\b(u),
\end{eqnarray*}
for any $a,b\in A$  and $u,v\in M$.  We denote it by $A\ltimes_{l_A,r_A} M$.
\end{proposition}

 {\bf Proof.}  Similar to Proposition 3.1.    \hfill $\square$

 \begin{lemma}
 Let $(M, l_A, r_A,  \b)$ be a module of a Hom-left-symmetric conformal algebra $(A, \a)$. Then

 (1) $l_A: A\rightarrow Cend(M)$ is a representation of the sub-adjacent Hom-Lie conformal algebra $(\mathfrak{g}(A), \a)$.

 (2) $\rho=l_A-r_A$ is a representation of the  Hom-Lie conformal algebra $(\mathfrak{g}(A), \a)$.

 (3) For any representation $\sigma: \mathfrak{g}(A)\rightarrow Cend(M)$ of the  Hom-Lie conformal algebra $(\mathfrak{g}(A), \a)$, $(M, \sigma, 0)$ is an $A$-module.
\end{lemma}

{\bf Proof.}  We only show that  (2) holds. To prove (2), we need to show that $\rho([a_\lambda b])_{\lambda+\mu}=[\rho(a)_\lambda \rho(b)]$.In fact
\begin{eqnarray*}
&&[(l_A-r_A)(a)_{\lambda}, (l_A-r_A)(b)_{\mu}]\b(v)- (l_A-r_A)([a_\lambda b])_{\lambda+\mu}\b(v)\\
&=& (l_A-r_A)(\a(a))_{\lambda}(b_\mu v-v_{-\mu-\partial}b)-(l_A-r_A)(\a(b))_{\mu}(a_\lambda v-v_{-\mu-\partial}a)\\
&&-(l_A-r_A)(a_\lambda b-b_{-\mu-\lambda}a)\b(v)\\
&=& \a(a)_{\lambda}(b_\mu v)- (b_\mu v)_{-\lambda-\partial}\a(a)-    \a(a)_{\lambda}(v_{-\mu-\partial}b)+(v_{-\mu-\partial}b)_{-\lambda-\partial} \a(a)\\
&& -\a(b)_{\mu}(a_\lambda v)+ (a_\lambda v)_{-\mu-\partial}\a(b)+  \a(b)_{\mu}(v_{-\lambda-\partial}a)-(v_{-\lambda-\partial}a)_{-\mu-\partial} \a(b)\\
&& -(a_\lambda b)_{\lambda+\mu} \b(v)+\b(v)_{-\lambda-\mu-\partial}(a_\lambda b)+(b_{-\lambda-\partial}a)_{\lambda+\mu}\b(v)-\b(v)_{-\lambda-\mu-\partial}(b_{-\lambda-\partial} a)\\
&=& (\a(a)_{\lambda}(b_\mu v)-\a(b)_{\mu}(a_\lambda v)-(a_\lambda b)_{\lambda+\mu} \b(v)+(b_{-\lambda-\partial}a)_{\lambda+\mu}\b(v))\\
&& (-(b_\mu v)_{-\lambda-\partial}\a(a)+(v_{-\mu-\partial}b)_{-\lambda-\partial} \a(a)+ \a(b)_{\mu}(v_{-\lambda-\partial}a)-\b(v)_{-\lambda-\mu-\partial}(b_{-\lambda-\partial} a))\\
&& (-\a(a)_{\lambda}(v_{-\mu-\partial}b)+(a_\lambda v)_{-\mu-\partial}\a(b)-(v_{-\lambda-\partial}a)_{-\mu-\partial} \a(b)+\b(v)_{-\lambda-\mu-\partial}(a_\lambda b))\\
&=& 0.
\end{eqnarray*}
And finish this proof.   \hfill $\square$

\begin{proposition}
Let $(A, \a)$ be a Hom-left-symmetric conformal algebra and  $(M, l_A, r_A,  \b)$ be a  finite module. Let $l^{\ast}_A,r^{\ast}_A: A\rightarrow Cend(M^{\ast c})$ be two $\mathbb{C}[\partial]$-module homomorphisms given by
\begin{eqnarray}
&& \a^{\ast}(a^{\ast})(b)=a^{\ast}(\a(y)), \b^{\ast}(u^{\ast})(v)=u^{\ast}(\b(v)),\\
&& (l^{\ast}_A(a)_{\lambda}f)_{\mu} u=-f_{\mu-\lambda}(l_A(a)_{\lambda}u), (r^{\ast}_A(a)_{\lambda}f)_{\mu}u=-f_{\mu-\lambda}(r_A(a)_{\lambda}u).
\end{eqnarray}
If in addition,
 \begin{eqnarray}
&&\b(l_A(\a(a))m)=l_A(a)\b(m), \b(r_A(\a(a))m)=r_A(a)\b(m),\\
&& \b(l_A(a_\lambda b)_{\lambda+\mu}v)-l_A(a)_\lambda(l_A(\a(b))v)=\b(l_A(b_{\mu}a)_{\lambda+\mu}v)-l_A(b)_{\mu}  (l_A(\a(a))_{\lambda}v),~~~\\
&& r_A(b)_{-\lambda-\mu-\partial}(l_A(\a(a))_\lambda v)-l_A(a)_\lambda(r_A(\a(b))_{-\mu-\partial}v)\nonumber\\
&& =r_A(b)_{-\lambda-\partial\mu}(r_A(\a(a))_{\lambda}v)-\b(r_A(a_\lambda b)_{-\mu-\partial} v),
 \end{eqnarray}
for any $a, b\in A, a^{\ast}\in A^{\ast},  f\in M^{\ast c}$ and $u\in M$.  Then $(M^{\ast c}, l^{\ast}_A-r^{\ast}_A, -r^{\ast}_A, \b^{\ast})$ is an $A$-module.
\end{proposition}

{\bf Proof. }  Let $\rho^{\ast}=l^{\ast}_A-r^{\ast}_A$. According to (3.6) and (3.7), we have
     \begin{eqnarray*}
    && \b^{\ast}(\rho^{\ast}(a)_{\lambda}f)_{\mu} u\\
     &=& (l^{\ast}_A(a)_{\lambda}f)_{\mu} \b(u)-(r^{\ast}_A(a)_{\lambda}f)_{\mu} \b(u)\\
     &=& -f_{\mu-\lambda}(l_A(a)_{\lambda}\b(u))+f_{\mu-\lambda}(r_A(a)_{\lambda}\b(u))\\
     &=& f_{\mu-\lambda}((r_A(a)-l_A(a))_{\lambda}\b(u)),
     \end{eqnarray*}
 and
 \begin{eqnarray*}
&& \rho^{\ast} (\a(a))_{\lambda}\b^{\ast}(f)_{\mu}u\\
&=& l^{\ast}_A(\a(a))_{\lambda}\b^{\ast}(f)_{\mu}u- r^{\ast}_A(\a(a))_{\lambda}\b^{\ast}(f)_{\mu}u\\
&=& f_{\mu-\lambda}(\b(r_A(\a(a))-l_A(\a(a)))_{\lambda}v),
 \end{eqnarray*}
 that is $(3.1)$ holds for   $l^{\ast}_A-r^{\ast}_A$. Similarly, (3.1) holds for $-r^{\ast}_A$. According to (3.8)-(3.10), we have
 \begin{eqnarray*}
&& -(\rho^{\ast}_A(\a(a))_{\lambda}(r^{\ast}_A(b)_{-\mu-\partial}f))_{\nu}u +r^{\ast}_A(\a(b))_{-\lambda-\mu-\partial}(\rho^{\ast}_A(a)_{\lambda}f))_{\nu}u\\
&&+(r^{\ast}_A(\a(b))_{-\lambda-\mu-\partial}(r^{\ast}_A(a)_{\lambda}f)_{\nu}u+(r^{\ast}_A(a_{\lambda}b)_{-\nu-\partial}\b^{\ast}(f))_{\nu}u\\
&=& (\rho^{\ast}_A(\a(a))_{\lambda}f)_{\lambda+\mu} (r_A(b)_{-\lambda-\mu+\nu}u)-(r^{\ast}_A(\a(b))_{-\mu-\partial}f)_{\nu-\lambda}(\rho_A(a)_{\lambda}(v))\\
&& -(r_A^{\ast}(a)_{\lambda}f)_{\lambda+\mu}(r_A(\a(b))_{-\lambda-\mu+\nu}u)-\b^{\ast}(f_{\mu}(r_A(a_{\lambda}b)_{\nu-\mu}u))\\
&=& -f_{\mu}(r_A(b)_{\nu-\lambda-\mu}((l_A(\a(a))-r_A(\a(a)))_{\lambda}u))+f_{\mu}((l_A(a)-r_A(a))_{\lambda}(r_A(\a(b))_{-\lambda-\mu+\nu}u))\\
&& +f_{\mu}(r_A(a)_{\lambda}(r_A(\a(b))_{-\lambda-\mu+\nu}u))-f_{\mu}(\b(r_A(a_\lambda b)_{\nu-\mu}u))\\
&=& -f_{\mu}(r_A(b)_{\nu-\lambda-\mu}((l_A(\a(a))-r_A(\a(a)))_{\lambda}u))+f_{\mu}((l_A(a))_{\lambda}(r_A(\a(b))_{-\lambda-\mu+\nu}u))\\
&& -f_{\mu}(\b(r_A(a_\lambda b)_{\nu-\mu}u))=0.
 \end{eqnarray*}
 Similarly, by $[\rho^{\ast}_A(a)_{\lambda}, \rho^{\ast}_A(b)_{\mu}]=\rho^{\ast}([a_{\lambda}b]), \forall a,b\in A$, (3.2) holds. Therefore, we have that $(M^{\ast c}, l^{\ast}_A-r^{\ast}_A, -r^{\ast}_A, \b^{\ast})$ is an $A$-module.    \hfill $\square$

\begin{proposition}

Let $(A, \a)$  and $(B,  \g )$ be two  Hom-left-symmetric conformal algebras. Suppose that there are $\mathbb{C}[\partial]$-module homomorphisms $l_A,r_A: A\rightarrow gc(B)$ and $l_B,r_B:B\rightarrow gc(A)$ such that $(B,l_A,r_A, \g)$ is an $A$-module and $(A,l_B,r_B, \a)$ is a $B$-module satisfying the following conditons
\begin{eqnarray}
&& r_A(\a(x))_{-\lambda-\mu-\partial}(a_\lambda b)=r_A(l_B(b)_{\mu}x)_{-\lambda-\partial}\g(a)-r_A(l_B(a)_{\lambda}x)_{-\mu-\partial}\g(b)\nonumber\\
&& +\g(a)_{\lambda}(r_A(x)_{-\mu-\partial}b)-\g(b)_{\mu}(r_A(x)_{-\lambda-\partial}a),\\
&& l_A(\a(x))_{\lambda}(a_\mu b)=-l_A(l_B(a)_{\mu}x-r_B(a)_{-\lambda-\partial}x)_{\lambda+\mu}\g(b)+(l_A(x)_{\lambda}a\nonumber\\
&&-r_A(x)_{-\mu-\partial}a)_{\lambda+\mu}\g(b)+r_A(r_B(b)_{-\lambda-\partial}x)_{-\mu-\partial}\g(a)+\g(a)_{\mu}(l_A(x)_{\lambda}b),\\
&& r_B(\g(a))_{-\lambda-\mu-\partial}(x_\lambda y)=r_B(l_A(y)_{\mu}a)_{-\lambda-\partial}\a(x)-r_B(l_A(x)_{\lambda}a)_{-\mu-\partial}\a(y)\nonumber\\
&& +\a(x)_{\lambda}(r_B(a)_{-\mu-\partial}y)-\a(y)_{\mu}(r_B(a)_{-\lambda-\partial}x),\\
&& l_B(\g(a))_{\lambda}(x_\mu y)=-l_B(l_A(x)_{\mu}a-r_A(x)_{-\lambda-\partial}a)_{\lambda+\mu}\a(y)+(l_B(a)_{\lambda}x\nonumber\\
&&-r_B(a)_{-\mu-\partial}x)_{\lambda+\mu}\a(y)+r_B(r_A(y)_{-\lambda-\partial}a)_{-\mu-\partial}\a(x)+\a(x)_{\mu}(l_B(a)_{\lambda}y),
\end{eqnarray}
for any $x,y\in A$ and $a,b\in B$. Then there is a Hom-left-symmetric conformal algebra on the $\mathbb{C}[\partial]$-module $A\oplus B$ given by
\begin{eqnarray*}
&&(x+a)_{\lambda}(y+b)=(x_\lambda y+l_B(a)_{\lambda}y+r_{B}(b)_{-\lambda-\partial}x)+ (a_\lambda b+l_A(x)_{\lambda}b+r_A(y)_{-\lambda-\partial}a),\\
&&(\a+\g)(x+a)=\a(x)+\g(a),
\end{eqnarray*}
for any $x,y\in A$ and $a,b\in B$.  We denote this Hom-left-symmetric conformal algebra by $(A\bowtie B, \a+\g)$. And $(A, B,l_A,r_A, l_B,r_B, \a,\g)$ satisfying the above conditions is called  a matched pair of Hom-left-symmetric conformal algebras.
\end{proposition}

{\bf Proof.}  Similar to Proposition 3.2.    \hfill $\square$

\begin{corollary}
Let $(A, \a)$  and $(B,  \g )$ be two  Hom-left-symmetric conformal algebras. Suppose that there are $\mathbb{C}[\partial]$-module homomorphisms $\rho: A\rightarrow gc(B)$ and $\sigma:B\rightarrow gc(A)$ such that
\begin{eqnarray*}
&& \rho(\a(x))_{\lambda}(a\mu b)=(\rho(x)_{\lambda}a)_{\lambda+\mu} \g(b)-(\rho(\sigma(a))_{\mu}x)_{\lambda+\mu}\g(b)+\g(a)_{\mu}(\rho(x)_{\lambda}b),\\
&&\sigma(\g(a))_{\lambda}(x\mu y)=(\rho(a)_{\lambda}x)_{\lambda+\mu} \a(y)-(\sigma(\rho(x))_{\mu}a)_{\lambda+\mu}\a(y)+\a(x)_{\mu}(\sigma(a)_{\lambda}y),
\end{eqnarray*}
for any $x,y\in A$ and $a,b\in B$. Then the $\mathbb{C}[\partial]$-module  $A\oplus B$ given by
\begin{eqnarray*}
&&(x+a)_{\lambda}(y+b)=(x_\lambda y)+\sigma(a)_{\lambda}y)+ (a_\lambda b+\rho(x)_{\lambda}b),\\
&& (\a+\g)(x+a)=\a(x)+\g(a).
\end{eqnarray*}
Therefore, its sub-adjacent Hom-Lie conformal algebra is given by
\begin{eqnarray*}
&&[(x+a)_{\lambda}(y+b)]=(x_\lambda y-y_{-\lambda-\partial}x+\sigma(a)_{\lambda}y-\sigma(b)_{-\lambda-\partial}x)+ (a_\lambda b-b_{-\lambda-\partial}a\\
&&+\rho(x)_{\lambda}b-\rho(y)_{-\lambda-\partial}a),\\
&& (\a+\g)(x+a)=\a(x)+\g(a).
\end{eqnarray*}
for any $x,y\in A$ and $a,b\in B$.
\end{corollary}

\begin{corollary}
Let $(A, B,l_A,r_A, l_B,r_B, \a,\g)$ be  a matched pair of Hom-left-symmetric conformal algebras. Then $(\mathfrak{g}(A),\mathfrak{g}(B), l_A-r_A,l_B-r_B, \a,\g )$ is a  matched pair of Hom-lie conformal algebras.
\end{corollary}

\begin{theorem}
Let $(A, \a)$ be a Hom-left-symmetric conformal algebra which is free and finite as a $\mathbb{C}[\partial]$-module. Suppose there is another Hom-left-symmetric conformal algebra structure on the $\mathbb{C}[\partial]$-module $(A^{\ast c}, \a^{\ast})$. Then $(\mathfrak{g}(A), \mathfrak{g}(A^{\ast c}), L^{\ast}_A, L^{\ast}_{A^{\ast c}}, \a, \a^{\ast})$ is a matched pair of Hom-Lie conformal algebras if and only if $(A, A^{\ast c}, ad^{\ast}_A, -R^{\ast}_A, ad^{\ast}_{A^{\ast c}} -R^{\ast}_{A^{\ast c}}, \a, \a^{\ast})$ is a matched pair of Hom-left-symmetric conformal algebras.
\end{theorem}

{\bf Proof.}   $\Rightarrow$  let $l_A=ad^{\ast}_A, r_A=-R^{\ast}_A, l_B=l_{A^{\ast c}}=ad^{\ast}_{A^{\ast c}}, r_B=r_{A^{\ast c}}=-R^{\ast}_{A^{\ast c}}, \rho=L^{\ast}_A$ and $\sigma=L^{\ast}_{A^{\ast c}}$.  Let ${e_1, ..., e_n}$ be a $\mathbb{C}[\partial]$-basis of $A$, and ${e^{\ast}_1, ..., e^{\ast}_n}$ be a dual $\mathbb{C}[\partial]$-basis of $A^{\ast c}$ in the sense that $(e^{\ast}_j)_{\lambda}e_i=\delta_{ij}$. Set $e_{i\lambda}e_j=\sum^{n}_{k=1} P_k^{ij}(\lambda, \partial)e_k$ and $e^{\ast}_{i\lambda}e^{\ast}_j=\sum^{n}_{k=1} R_k^{ij}(\lambda, \partial)e^{\ast}_k$, where $P_k^{ij}(\lambda, \partial), R_k^{ij}(\lambda, \partial)\in \mathbb{C}[\lambda,\partial]$. Since
\begin{eqnarray*}
(L^{\ast}_A(e_i)_{\lambda}e^{\ast}_j)_\mu e_k=-e^{\ast}_{j \mu-\lambda}(e_{i \lambda}e_k)=-P_j^{ik}(\lambda, \mu-\lambda),
\end{eqnarray*}
we have
\begin{eqnarray*}
 L^{\ast}_A(e_i)_{\lambda}e^{\ast}_j=-\sum^{n}_{k=1} P_j^{ik}(\lambda, \partial)e_k^{\ast}.
\end{eqnarray*}
Similarly, we have
\begin{eqnarray*}
&&R^{\ast}_A(e_i)_{\lambda}e^{\ast}_j=-\sum^{n}_{k=1} P_j^{ki}(\partial, -\lambda-\partial)e_k^{\ast},\\
&& L^{\ast}_{A^{\ast c}}(e^{\ast}_i)_{\lambda}e_j=-\sum^{n}_{k=1} R_j^{ik}(\lambda, -\lambda-\partial)e_k,\\
&& R^{\ast}_{A^{\ast c}}(e^{\ast}_i)_{\lambda}e_j^{\ast}=-\sum^{n}_{k=1} R_j^{ki}(\partial, -\lambda-\partial)e_k.
\end{eqnarray*}
By direct computation , we have
 \begin{eqnarray*}
&&(L^{\ast}_A(e_i)_{\lambda}([e^{\ast}_{j\mu}e^{\ast}_k]))_{\nu}e_t=-(R^{\ast}_A(e_t)_{-\eta-\theta-\partial}([e^{\ast}_{j\eta}e^{\ast}_k]))_{\omega}e_i\\
 && -([(L^{\ast}_A(e_i)_{\lambda}e^{\ast}_j)_{\lambda+\mu}e^{\ast}_k])_{\nu}e_t=(R^{\ast}_A(ad^{\ast}_{A^{\ast c}}(e^{\ast}_k)_{\theta}e_t)_{-\eta-\partial}e^{\ast}_j)_{\omega}e_i\\
 && ([e^{\ast}_{j\mu}(L^{\ast}_A(e_i)_{\lambda}e^{\ast}_k)])_{\nu}e_{t}=(R^{\ast}_A(ad^{\ast}_{A^{\ast c}}(e^{\ast}_j)_{\eta}e_t)_{-\theta-\partial}e^{\ast}_k)_{\omega}e_i\\
 && (L^{\ast}_A(L^{\ast}_{A^{\ast c}}(e^{\ast}_j)_{-\lambda-\partial}e_i)_{\lambda+\mu}e^{\ast}_k)_{\nu}e_t=(e^{\ast}_{j\eta}(R^{\ast}_A(e_t)_{-\theta-\partial}e^{\ast}_k))_{\omega}e_i\\
 && (L^{\ast}_A(L^{\ast}_{A^{\ast c}}(e^{\ast}_k)_{-\lambda-\partial} e_i)_{-\mu-\partial}e^{\ast}_j)_{\nu}e_t=(e^{\ast}_{k \theta}(R^{\ast}_A(R^{\ast}_{A}(e_t)_{-\eta-\partial}e^{\ast}_j)))_{\omega}e_i,
 \end{eqnarray*}
by letting $\eta=\mu, \omega=-\lambda$ and $\theta=\nu-\lambda-\mu$. Therefore, (3.1)$\Leftrightarrow$ (3.12), Similarly, we get
\begin{eqnarray*}
(3.1) \Leftrightarrow  (3.12) \Leftrightarrow (3.15)\\
(3.2) \Leftrightarrow  (3.13) \Leftrightarrow (3.14)
\end{eqnarray*}
$\Leftarrow$ It follows from Corollary 3.10.  \hfill $\square$

\section{ Hom-left-symmetric conformal  bialgebras}
\def\theequation{\arabic{section}. \arabic{equation}}
\setcounter{equation} {0}
In this section, we introduce the notion of  a Hom-left-symmetric conformal  bialgebra, which is equivalent to the parak\"{a}hler Hom-Lie conformal algebra. In particular, we investigate the coboundary Hom-left-symmetric conformal  bialgebras.

\begin{definition}
A Hom-left-symmetric conformal  coalgebra $(A, \a)$ is a $\mathbb{C}[\partial]$-module equipped with an even linear endomorphism
$\a$ such that $\a\partial=\partial\a$, and  a $\mathbb{C}$-module homomorphism $\Delta: A\rightarrow A\o A$ such that
\begin{eqnarray*}
(\a\o \Delta)\Delta(x)-\tau_{12}(\a\o \Delta)\Delta(x)=( \Delta\o \a)\Delta(x)-\tau_{12}(\Delta \o\a)\Delta(x),
\end{eqnarray*}
where $\tau_{12}(x\o y\o z)=y\o x\o z$ for any $x,y,z\in A$.
\end{definition}

\begin{proposition}
Let $(A, \a, \Delta)$ be a finite Hom-left-symmetric conformal  coalgebra. Then $A^{\ast c}=Chom(A, \mathbb{C})$ is a Hom-left-symmetric conformal  algebra endowed with the following product
\begin{eqnarray*}
(f_\mu g)_{\lambda}(r)=\sum f_\mu(r_{(1)})g_{\lambda-\mu}(r_{(2)})=(f\o g)_{\mu, \lambda-\mu}(\Delta(r)),
\end{eqnarray*}
where $\Delta(r)=\sum r_{(1)}\o r_{(2)}$.
\end{proposition}

{\bf Proof.} It is easy to check that (2.1) holds, we only check that
\begin{eqnarray}
(f_{\lambda}g)_{\lambda+\mu}\a^{\ast}(h)-\a^{\ast}(f)_{\lambda}(g_{\mu}h)=(g_{\mu}f)_{\lambda+\mu} \a^{\ast}(h)- \a^{\ast}(g)_{\mu}(f_\lambda h).
\end{eqnarray}
For this, we calculate
\begin{eqnarray*}
&&((f_{\lambda}g)_{\lambda+\mu}\a^{\ast}(h)-\a^{\ast}(f)_{\lambda}(g_{\mu}h))_{\nu}(r)\\
&=& \sum (f_{\lambda}g)_{\lambda+\mu}(r_{(1)}) \a^{\ast}(h)_{\nu-\lambda-\mu}(r_{(2)})-\sum \a^{\ast}(f)_{\lambda}(r_{(1)})(g_{\mu}h))_{\nu-\lambda}(r_{(2)})\\
&=& \sum f_{\lambda}(r_{(11)})g_{\mu}(r_{(12)})h_{\nu-\lambda-\mu}(\a(r_{(2)}))-\sum f_{\lambda}(\a(r_{(1)}))g_{\mu}(r_{(21)})h_{\nu-\lambda-\mu}(r_{(22)}),\\
&&((g_{\mu}f)_{\lambda+\mu} \a^{\ast}(h)- \a^{\ast}(g)_{\mu}(f_\lambda h))_{\nu}(r)\\
&=& \sum (g_{\mu}f)_{\lambda+\mu}(r_{(1)})  \a^{\ast}(h)_{\nu-\lambda-\mu}(r_{(2)})-\sum \a^{\ast}(g)_{\mu}(r_{(1)})(f_{\lambda}h))_{\nu-\mu}(r_{(2)})\\
&=& \sum f_{\lambda}(r_{(12)})g_{\mu}(r_{(11)})h_{\nu-\lambda-\mu}(\a(r_{(2)}))-\sum f_{\lambda}(r_{(2)})g_{\mu}(\a(r_{(1)}))h_{\nu-\lambda-\mu}(r_{(22)}),
\end{eqnarray*}
where $\delta(r_{(1)})=\sum r_{(11)}\o r_{(12)}$ and $\delta(r_{(2)})=\sum r_{(21)}\o r_{(22)}$. Then $A^{\ast c}=Chom(A, \mathbb{C})$ is a Hom-left-symmetric conformal  algebra.
 \hfill $\square$

\begin{proposition}
Let $(A, \a)$ be a  Hom-left-symmetric conformal  algebra which is free of finite rank, that is $A=\sum_{i=1}^{n}\mathbb{C}[\partial]e_i$. Then $A^{\ast c}=Chom(A, \mathbb{C})=\sum_{i=1}^{n}\mathbb{C}[\partial]e_i^{\ast}$, where ${e_i^{\ast}}$ is a dual $\mathbb{C}[\partial]$-basis of $A^{\ast c}$ in the sense that $(e^{\ast}_i)_{\lambda}e_j=\delta_{ij}$, is a Hom-left-symmetric conformal  coalgebra with the following coproduct
\begin{eqnarray*}
\delta(f)=\sum_{i,j}f_{\mu}(e_{i\lambda}e_j)(e_{i}^{\ast}\o e_{j}^{\ast})|_{\lambda=\partial\o 1, \mu=-\partial\o 1-1\o \partial}.
\end{eqnarray*}
\end{proposition}

{\bf Proof.} Similar to \cite{Hong2015}.  \hfill $\square$

Let $(A, \a)$ be a  Hom-left-symmetric conformal  algebra, we know that $(A\o A, \phi_A)$ is a module over $\mathfrak{g}(A)$, we consider the reduced complex for the $\mathfrak{g}(A)$-module $(A\o A, \phi_A)$. The details about cohomology of Hom-conformal algebras can be referred to \cite{Zhao2018}. A 1-cochain is the map $\delta: A\rightarrow (A\o A)[\lambda]$ such that $\delta_{\lambda}(\partial a)=-\lambda\delta_{\lambda}( a)$. 
A 1-cochain is the map $\delta: A\rightarrow A\o A$ such that $\delta$ is a $\mathbb{C}[\partial]$-module homomorphism. The condition $d\delta=0$ as
\begin{eqnarray}
\delta(\a([a_{\lambda}b]))=\phi_A(a)_{\lambda}\delta(b)-\phi_A(b)_{-\lambda-\partial}\delta(a), a,b\in A.
\end{eqnarray}
If $\delta$ is a $\mathbb{C}[\partial]$-module homomorphism and satisfies (4.2), we say $\delta$ is a 1-cocycle of $\mathfrak{g}(A)$ associated to $(A\o A, \phi_A)$.

\begin{theorem}
Let  $(A, \a)$ be a  Hom-left-symmetric conformal  algebra which is free and finite as a $\mathbb{C}[\partial]$-module, whose $\lambda$-product is obtained from a $\mathbb{C}[\partial]$-module homomorphism $\psi: A^{\ast c}\rightarrow A^{\ast c}\o A^{\ast c}$.  Suppose there is another Hom-left-symmetric conformal algebra structure on the $\mathbb{C}[\partial]$-module $(A^{\ast c}, \a^{\ast})$ obtained from a $\mathbb{C}[\partial]$-module homomorphism $\varphi: A\rightarrow A\o A$. Then $(\mathfrak{g}(A), \mathfrak{g}(A^{\ast c}), L^{\ast}_A, L^{\ast}_{A^{\ast c}}, \a, \a^{\ast})$ is a matched pair of Hom-Lie conformal algebras if and only if $\varphi: A\rightarrow A\o A$ is a 1-cocycle of $\mathfrak{g}(A)$ associated to $L_A\o \a+\a\o ad_A$  and $\psi: A^{\ast c}\rightarrow A^{\ast c}\o A^{\ast c}$ is a 1-cocycle of $\mathfrak{g}(A^{\ast c})$ associated to $L_{A^{\ast c}}\o \a^{\ast}+\a^{\ast}\o ad_{A^{\ast c}}$.
\end{theorem}
{\bf Proof.}    Let ${e_1, ..., e_n}$ be a $\mathbb{C}[\partial]$-basis of $A$ and ${e^{\ast}_1, ..., e^{\ast}_n}$ the dual $\mathbb{C}[\partial]$-basis of $A^{\ast c}$ in the sense that $(e^{\ast}_j)_{\lambda}e_i=\delta_{ij}$. Set $e_{i\lambda}e_j=\sum^{n}_{k=1} P_k^{ij}(\lambda, \partial)e_k$ and $e^{\ast}_{i\lambda}e^{\ast}_j=\sum^{n}_{k=1} R_k^{ij}(\lambda, \partial)e^{\ast}_k$, where $P_k^{ij}(\lambda, \partial), R_k^{ij}(\lambda, \partial)\in \mathbb{C}[\lambda,\partial]$. By Proposition 4.3, we have
\begin{eqnarray*}
&& \varphi (e_k)=\sum_{i,j}Q_k^{ij}(\partial\o 1, 1\o \partial) e_i\o e_j, ~~\mbox{where}~~ Q_k^{ij}(x, y)=R_k^{ij}(x, -x-y),\\
&& \psi(e_k^{\ast})=\sum_{i,j}S_k^{ij}(\partial\o 1, 1\o \partial) e^{\ast}_i\o e^{\ast}_j, ~\mbox{where}~~ S_k^{ij}(x, y)=P_k^{ij}(x, -x-y),
\end{eqnarray*}
Furthermore, we have
\begin{eqnarray*}
&&\varphi([e_{i\lambda}e_j])=[L_A(e_i)_{\lambda}\o \a+\a\o ad_A(e_i)_{\lambda}]\varphi(e_j)\\
&&-[L_A(e_j)_{\lambda}\o \a+\a\o ad_A(e_j)_{\lambda}]\varphi(e_i).
\end{eqnarray*}
Let
$
\a(e_i)=\sum_{s}f(\partial)e_s, \a(e_t)=\sum_{n}g(\partial)e_n,  \a(e_l)=\sum_{n}h(\partial)e_m.
$
Hence the coefficient of $e_{m}\o e_n$ gives the following relation
\begin{eqnarray*}
&&\sum _{k} (p_{k}^{ij} (\lambda, \partial^{\otimes^2})-P_{k}^{ji}(-\lambda-\partial^{\otimes^2}, \partial^{\otimes^2})) R_{k}^{st}(\partial\o 1, -\partial^{\otimes^2})f(\partial)\\
&=& \sum_{k} (p_{s}^{ik} (\lambda, \partial\o 1)R_{j}^{kt}(\lambda+\partial\o 1,-\lambda-\partial^{\otimes^2})g(\partial)\\
&& +\sum_{k} h(\partial) (P_{t}^{ik}(\lambda, 1\o \partial)-P_{t}^{ki}(-\lambda-1\o \partial, 1\o \partial)) R_{j}^{sk}(1\o \partial, -\lambda-\partial^{\otimes^2})
\end{eqnarray*}
\begin{eqnarray*}
&& -\sum_{k}(p_{s}^{jk} (-\lambda-\partial^{\otimes^2}, \partial\o 1)R_{i}^{kt}(-\lambda-1\o \partial,\lambda)g(\partial)\\
&& -\sum_{k} h(\partial) (P_{t}^{jk}(-\lambda-\partial^{\otimes^2}, 1\o \partial)-P_{t}^{kj}(\lambda+\partial\o1, 1\o \partial)) R_{i}^{sk}( \partial\o 1, \lambda),
\end{eqnarray*}
where $\partial^{\otimes^2}=\partial\o 1+1+\partial$, which is precisely the relation given by the coefficient of $e_n$ in
\begin{eqnarray*}
&&-L^{\ast}_{A^{\ast c}}(\a^{\ast}(e^{\ast}_m))_{-\lambda-\mu-\partial}[e_{i\lambda}e_j]=L^{\ast}_{A^{\ast c}}(L^{\ast}_A(e_i)_{\lambda}e^{\ast}_m)_{-\mu-\partial}\a(e_j) -[\a(e_i)_{\lambda}(L^{\ast}_{A^{\ast c}}(e^{\ast}_m))_{-\mu-\partial}e_j]\\
&&-L^{\ast}_{A^{\ast c}}(L^{\ast}_A(e_j)_{\mu}e^{\ast}_m)_{-\lambda-\partial}\a(e_i)+[\a(e_j)_{\mu}(L^{\ast}_{A^{\ast c}}(e^{\ast}_m))_{-\lambda-\partial}e_i].
\end{eqnarray*}
 It is just the condition for  $(\mathfrak{g}(A), \mathfrak{g}(A^{\ast c}), L^{\ast}_A, L^{\ast}_{A^{\ast c}}, \a, \a^{\ast})$  is a matched pair of Hom-Lie conformal algebras.  Then $\varphi: A\rightarrow A\o A$ is a 1-cocycle of $\mathfrak{g}(A)$ if and only if (3.2) holds in the case $\rho=L^{\ast}_A$ and $\sigma=L^{\ast}_{A^{\ast c}}$. Similarly, $\psi: A^{\ast c}\rightarrow A^{\ast c}\o A^{\ast c}$ is a 1-cocycle of $\mathfrak{g}(A^{\ast c})$ if and only if (3.1) holds when $\rho=L^{\ast}_A$ and $\sigma=L^{\ast}_{A^{\ast c}}$.  \hfill $\square$

\begin{definition}
Let $(A,\a)$ be a finitely generated  $\mathbb{C}[\partial]$-module. A  Hom-left-symmetric conformal  bialgebra on $(A, \a)$ is a six-tuple $(A, A^{\ast c}, \a, \a^{\ast}, \varphi, \psi)$, where $\varphi: A\rightarrow A\o A$ and $\psi: A^{\ast c}\rightarrow A^{\ast c}\o A^{\ast c}$ are two $\mathbb{C}[\partial]$-module homomorphisms, $(A, \a, \varphi)$ and $(A^{\ast}, \a^{\ast}, \psi)$ are two Hom-left-symmetric conformal coalgebras, satisfying the following conditions:

(1) $\varphi$ is a 1-cocycle of $\mathfrak{g}(A)$ associated to $L_{A}\o \a+\a\o ad_{A}$,

(2) $\psi$ is a 1-cocycle of $\mathfrak{g}(A^{\ast c})$ associated to $L_{A^{\ast c}}\o \a^{\ast}+\a^{\ast}\o ad_{A^{\ast c}}$.
\end{definition}

We denote this Hom-left-symmetric conformal  bialgebra by $(A, A^{\ast c}, \varphi, \psi)$.

By Theorem 3.11 and Theorem 4.4, we have the following proposition.
\begin{proposition}
Let $(A, \a)$ be a  Hom-left-symmetric conformal  algebra which is free and finite as a $\mathbb{C}[\partial]$-module. Suppose that there is also a Hom-left-symmetric conformal  algebra structure on $(A^{\ast c}, \a^{\ast})$. Then the following conditions are equivalent:

(1) $(A, A^{\ast c}, \varphi, \psi, \phi)$ is a Hom-left-symmetric conformal  bialgebra on $(A, \a)$, where $\varphi, \psi$ are two $\mathbb{C}[\partial]$-module homomorphisms which are obtained from the $\lambda$-products of $(A^{\ast c}, \a^{\ast})$ and $(A, \a)$, respectively.

(2) $(\mathfrak{g}(A), \mathfrak{g}(A^{\ast c}), L^{\ast}_A, L^{\ast}_{A^{\ast c}}, \a, \a^{\ast})$ is a matched pair of Hom-Lie conformal algebras.

(3) $(A, A^{\ast c}, ad^{\ast}_A, -R^{\ast}_A, ad^{\ast}_{A^{\ast c}} -R^{\ast}_{A^{\ast c}}, \a, \a^{\ast})$ is a matched pair of Hom-left-symmetric conformal algebras.

(4) $(\mathfrak{g}(A)\bowtie \mathfrak{g}(A^{\ast c}, \mathfrak{g}(A), \mathfrak{g}(A^{\ast c}, \omega_{\lambda})$ is a parak\"{a}hler   Hom-Lie conformal algebra.
\end{proposition}

\begin{definition}
A Hom-left-symmetric conformal  bialgebra $(A, A^{\ast c}, \varphi, \psi)$ is called coboundary if $\varphi$ is a 1-coboundary of associated to $L_{A}\o \a+\a\o ad_{A}$, that is there exists a $r\in A\o A$ such that
\begin{eqnarray}
\varphi(a)=(L_A(a)\o \a+\a\o ad_A(a))_{\lambda}r|_{\lambda=-\partial^{\otimes^2}},
\end{eqnarray}
 for any $a\in A$, where $\partial^{\otimes^2}=\partial\o 1+1+\partial$.
\end{definition}
Let $r=\sum_i r_i\o l_i\in A\o A$. Set
\begin{eqnarray*}
r_{12}=\sum_i r_i\o l_i,~r_{21}=\sum_il_i\o r_i.
\end{eqnarray*}
Define
\begin{eqnarray*}
[[r,r]]&=&\sum_{i,j}(r_{i\mu}r_{j}\o \a(l_j)\o \a(l_i))|_{\mu=1\o 1\o \partial}-\sum_{i,j}(\a(l_j)\o r_{i\mu}r_{j}\o  \a(l_i))|_{\mu=1\o 1\o \partial}\\
&& -\sum_{i,j}(\a(r_j)\o  [l_{j\mu}r_{i}]\o \a(l_i))|_{\mu= \partial\o 1\o 1}+\sum_{i,j}([l_{j\mu}r_{i}]\o \a(r_j)\o \a(l_i))|_{\mu= 1\o \partial \o 1}\\
&& -\sum_{i,j}(\a(r_i)\o  \a(r_j)\o [l_{i\mu}l_{j}])|_{\mu= \partial\o 1\o 1}.
\end{eqnarray*}
Let $Q(x)_{\lambda}=L_A(x)_{\lambda}\o \a\o \a+\a\o L_A(x)_{\lambda}\o \a+\a\o \a\o ad_A(x)_{\lambda}$, $P(x)_{\lambda}=L_A(x)_{\lambda}\o \a+\a\o L_A(x)_{\lambda}$. Set
\begin{eqnarray*}
M(a)&=&\sum_{j}P(a_{\lambda}r_j)_{\mu}(r_{12}-r_{21})\o \a^{2}(l_j)|_{\lambda=-\partial^{\otimes^3}, \mu=-\partial^{\otimes^2}\o 1}\\
&& -\sum _{j} P(\a(a))_{\lambda}(P(r_j))_{\mu} (r_{12}-r_{21})\o \a^{2}(l_j)|_{\lambda=-\partial^{\otimes^3}, \mu=-1\o 1\o \partial},\\
J_{\delta}(a)&=&\a^{\otimes^3}(Q(a)_{\lambda}[[r,r]]_{\lambda=-\partial^{\otimes^3}})+M(a),
\end{eqnarray*}
where $\partial^{\otimes^3}=\partial\o 1\o 1+1\o \partial\o 1+ 1\o 1\o \partial$.
\begin{proposition}
Let $(A, \a)$ be a Hom-left-symmetric conformal  algebra and $r=\sum_i r_i\o l_i\in A\o A$. Assume that $\a^{\otimes^2}(r)=r$. Define a map $\delta:A \rightarrow A\o A$ by
\begin{eqnarray*}
 \delta(a)=(L_A(a)\o \a+\a\o ad_A(a))_{\lambda}r|_{\lambda=-\partial^{\otimes^2}}, 
\end{eqnarray*}
 for any $a\in A.$ Then $(A, \delta, \a)$ is a Hom-left-symmetric conformal coalgebra if and only if $J_{\delta}(a)=0$.
\end{proposition}

{\bf Proof.} For any $a\in A$, we  only need to prove that
\begin{eqnarray*}
 J_{\delta}(a)=(\a\o \delta)\delta(a)-\tau_{12}(\a\o \delta)\delta(a)-( \delta\o \a)\delta(a)-\tau_{12}( \delta\o \a)\delta(a).
\end{eqnarray*}
We will repeat using $\a^{\otimes^{2}}(r)=r$, by computation, we have
\begin{eqnarray*}
&&(\a\o \delta)\delta(a)\\
&=& (\a\o \delta)(\sum_{i}a_\lambda r_i\o \a(l_i)+\a(r_i)\o [a_\lambda l_i])|_{\lambda=-\partial^{\otimes^2}}\\
&=& (\a\o \delta)(\sum_{i} a_\lambda r_i\o \a(l_i)+\a(r_i)\o [a_\lambda l_i])\\
&=&\sum_{i}\a(a_\lambda r_i)\o \delta(\a(l_i))+\a^2(r_i)\o \delta([a_\lambda l_i])\\
&=& \sum_{i,j}\a(a)_\lambda r_i\o \delta(l_i)+\a^2(r_i)\o \delta([a_\lambda l_i])
\end{eqnarray*}
\begin{eqnarray*}
&=&  \sum_{i,j}  \a(a)_\lambda r_i \o l_{i\mu}r_j\o \a(l_j) + \a(a)_\lambda r_i \o  \a(r_j)\o [l_{i\mu}l_j] \\
&&+\a(r_i)\o [a_{\lambda}l_i]_{\mu}r_j\o \a(l_j)+\a(r_i)\o \a(r_j)\o [[a_{\lambda}l_i]_{\mu}l_j]|_{\lambda=-\partial^{\otimes^3}, \mu=-1\o \partial^{\otimes^2}}.
\end{eqnarray*}
Similarly, we can obtain
\begin{eqnarray*}
&&(\delta\o \a)\delta(a)\\
&=& \sum_{i,j}  (a_\lambda r_i)_{\mu}r_j \o \a(l_j)\o \a^2(l_i) + \a(r_j)\o [(a_\lambda r_i)_{\mu}r_j] \o  \a^2(l_i) \\
&&+\a(r_{i})_{\mu}r_j\o \a(l_j)\o  [a_{\lambda}l_i]+\a(r_j)\o [\a(r_{i})_{\mu}l_j]\o [a_\lambda l_i]|_{\lambda=-\partial^{\otimes^3}, \mu=-1\o \partial^{\otimes^2}}.
\end{eqnarray*}
Furthermore, we have
\begin{eqnarray*}
 J_{\delta}(a)=(C1)+(C2)+(C3),
\end{eqnarray*}
where
\begin{eqnarray*}
(C1)&=&\sum_{i,j}  ((a_\lambda r_j)_{\mu}r_i \o \a(l_i)\o \a^2(l_j) + \a(r_i)\o [(a_\lambda r_j)_{\mu}r_i] \o  \a^2(l_j))|_{\lambda=-\partial^{\otimes^3}, \mu=-1\o \partial^{\otimes^2}} \\
&& - \sum_{i,j}(\a(l_i)\o (a_\lambda r_j)_{\mu}r_i \o \a^2(l_j) + [(a_\lambda r_j)_{\mu}l_i]\o \a(r_i)\o   \a^2(l_j))|_{\lambda=-\partial^{\otimes^3}, \mu=-1\o \partial^{\otimes^2}}\\
&&-\sum_{i,j}  (\a(a_\lambda r_i) \o \a(l_{i\mu}r_j)\o \a^2(l_j) + \a^2(r_i)\o [a_{\lambda}l_i]_{\mu}\a(r_j)\o \a^2(l_j) )|_{\lambda=-\partial^{\otimes^3}, \mu=-1\o \partial^{\otimes^2}}\\
&& +\sum_{i,j}  (\a(l_{i\mu}r_j)\o \a(a_\lambda r_i) \o \a^2(l_j) + [a_{\lambda}l_i]_{\mu}\a(r_j)\o \a^2(r_i)\o\a^2(l_j)) |_{\lambda=-\partial^{\otimes^3}, \mu=-1\o \partial^{\otimes^2}}\\
(C2)&=&\sum_{i,j} \a(r_{j})_{\mu}r_i\o \a(l_i)\o  [a_{\lambda}l_j]+\a(r_i)\o [\a(r_{j})_{\mu}l_i]\o [a_\lambda l_j]|_{\lambda=-\partial^{\otimes^3}, \mu=-1\o \partial^{\otimes^2}}\\
&& \sum_{i,j} \a(l_i)\o \a(r_{j})_{\mu}r_i\o   [a_{\lambda}l_j]+[\a(r_{j})_{\mu}l_i]\o \a(r_i)\o  [a_\lambda l_j]|_{\lambda=-\partial^{\otimes^3}, \mu=-1\o \partial^{\otimes^2}},\\
(C3)&=&-\sum_{i,j} \a(a)_\lambda r_i \o  \a(r_j)\o [l_{i\mu}l_j] + \a(r_i)\o \a(r_j)\o [[a_{\lambda}l_i]_{\mu}l_j]|_{\lambda=-\partial^{\otimes^3}, \mu=-1\o \partial^{\otimes^2}},\\
&& +\sum_{i,j}  \a(r_j)\o \a(a)_\lambda r_i \o [l_{i\mu}l_j] + \a(r_j)\o \a(r_i)\o [[a_{\lambda}l_i]_{\mu}l_j] |_{\lambda=-\partial^{\otimes^3}, \mu=-1\o \partial^{\otimes^2}},
\end{eqnarray*}
On the other hand, we have
\begin{eqnarray*}
 &&\a^{\otimes^3}(Q(a)_{\lambda}[[r,r]]_{\lambda=-\partial^{\otimes^3}})  \\
 &=& \sum_{i,j} ( \a(a)_{\lambda}(r_{i\mu}r_j)\o \a^{2}(l_j)\o \a^{2}(l_i)|_{\mu=1\o 1\o \partial} +  \a(r_{i\mu}r_j)\o \a(a_{\lambda} l_j)\o \a^{2}(l_i) |_{\mu=1\o 1\o \partial}\\
&& +\a(r_{i\mu}r_j)\o \a^{2}(l_j)\o [\a(a)_\lambda \a(l_i)]|_{\mu=-\partial^{\otimes^2}\o 1}-\a(a_{\lambda}l_j)\o \a(r_{i\mu}r_j)\o \a^{2}(l_i)  |_{\mu=1\o 1\o \partial}\\
&& -\a^{2}(l_j)\o \a(a)_\lambda(r_{i\mu}r_j)\o \a^{2}(l_i)|_{\mu=1\o 1\o \partial}-\a^{2}(l_j)\o \a(r_{i\mu}r_j)\o [\a(a)_\lambda \a(l_i)]|_{\mu=-\partial^{\otimes^2}\o 1}
\end{eqnarray*}
\begin{eqnarray*}
&& -\a(a_\lambda r_j)\o \a([l_{j\mu}r_i])\o \a^{2}(l_i)|_{\mu=-1\o \partial^{\otimes^2}}-\a^{2}(r_j)\o \a(a)_{\lambda}[l_{j\mu}r_i]\o \a^{2}(l_i)|_{\mu=\partial\o 1\o 1}\\
&& -\a^{2}(r_j)\o \a([l_{j\mu}r_i])\o[\a(a)_\lambda \a(l_i)]|_{\mu=\partial\o 1\o 1}+\a(a)_{\lambda}[l_{j\psi}r_i]\o \a^{2}(r_j)\o \a^{2}(l_i)|_{\mu=1\o \partial\o 1}\\
&& +\a([l_{j\psi}r_i])\o \a(a_{\lambda}r_j)\o \a^{2}(l_i)|_{\mu=-1\o 1\o \partial-\partial\o 1\o 1}+\a([l_{j\mu}r_i])\o \a^{2}(r_j)\o [\a(a)_{\lambda}\a(l_i)]|_{\mu=1\o \partial\o 1}\\
&& -\a(a_\lambda r_i)\o \a^{2}(r_j)\o\a([l_{i\mu}l_j])| _{\mu=-1\o \partial^{\otimes^2}}-\a^{2}(r_i)\o \a(a_{\lambda}r_j) \o \a([l_{i\mu}l_j])|_{\mu=\partial\o 1\o 1}\\
&& -\a^{2}(r_i)\o\a^{2}( r_j)\o [\a(a)_{\lambda}[l_{i\mu}l_j]]|_{\mu=\partial\o 1\o 1})|_{\lambda=-\partial^{\otimes^3}},\\
M(a)&=& \sum_{i,j} ((a_{\lambda}r_{j})_{\mu}r_i\o \a(l_i)\o \a^{2}(l_j)|_{\mu=\partial^{\otimes^2}\o 1}+\a(r_i)\o (a_\lambda r_j)_{\mu}l_i\o \a^{2}(l_j)|_{\mu=\partial^{\otimes^2}\o 1}\\
&& -(a_{\lambda}r_{j})_{\mu}l_i\o \a(r_i)\o \a^{2}(l_j)|_{\mu=\partial^{\otimes^2}\o 1}-\a(l_i)\o (a_\lambda r_j)_{\mu}r_i\o \a^{2}(l_j)|_{\mu=\partial^{\otimes^2}\o 1}\\
&& -\a(a)_{\lambda}(r_{j\mu}r_i)\o \a^{2}(l_i)\o \a^{2}(l_j)|_{\mu=1\o 1\o \partial}-\a(r_{j\mu}r_i)\o \a(a_{\lambda}l_i)\o \a^{2}(l_j)|_{\mu=1\o 1\o \partial}\\
&& -\a(a)_\lambda \a(r_i)\o \a(r_{j\mu}l_i)\o \a^{2}(l_j)|_{\mu=1\o 1\o \partial}-\a^{2}(r_i)\o \a(a)_{\lambda}(r_{j\mu}l_i)\o \a^{2}(l_j)|_{\mu=1\o 1\o \partial}\\
&&\a(a)_{\lambda}(r_{j\mu}l_i)\o \a^{2}(r_i)\o \a^{2}(l_j)|_{\mu=1\o 1\o \partial}+\a(r_{j\mu}l_i)\o \a(a)_{\lambda}\a(r_i)\o \a^{2}(l_j)|_{\mu=1\o 1\o \partial}\\
&& \a(a_\lambda l_i)\o \a(r_{j\mu}r_i)\o \a^{2}(l_j)|_{\mu=1\o 1\o \partial}+\a^{2}(l_i)\o \a(a)_{\lambda}(r_{j\mu}r_i)\o \a^{2}(l_j)|_{\mu=1\o 1\o \partial})|_{\lambda=-\partial^{\otimes^3}}.
\end{eqnarray*}
After rearranging the terms suitably, the sum of the terms whose third component is $\a^{2}(l_j)$ in $J_{\delta}(a)$ is
\begin{eqnarray*}
(D1)=(D11)+(D12)+(D13)+(D14)+(D15)+(D16)+(D17),
\end{eqnarray*}
where
\begin{eqnarray*}
(D11)&=& \sum_{i,j} \a(a)_{\lambda}(r_{i\mu}r_j)\o \a^{2}(l_j)\o \a^{2}(l_i)|_{\mu=1\o 1\o \partial}+(a_{\lambda}r_{j})_{\mu}r_i\o \a(l_i)\o \a^{2}(l_j)|_{\mu=\partial^{\otimes^2}\o 1}\\
&& -\a(a)_{\lambda}(r_{i\mu}r_j)\o \a^{2}(l_j)\o \a^{2}(l_i)|_{\mu=1\o 1\o \partial}|_{\lambda=-\partial^{\otimes^3}}\\
&=& \sum_{i,j}  (a_{\lambda}r_{j})_{\mu}r_i\o \a(l_i)\o \a^{2}(l_j)|_{\lambda=-\partial^{\otimes^3},\mu=\partial^{\otimes^2}\o 1}\\
(D12)&=& \sum_{i,j}-\a(l_i)\o (a_\lambda r_j)_{\mu}r_i\o \a^{2}(l_j)|_{\mu=\partial^{\otimes^2}\o 1}+\a^{2}(l_i)\o \a(a)_{\lambda}(r_{j\mu}r_i)\o \a^{2}(l_j)|_{\mu=1\o 1\o \partial}\\
&&-\a^{2}(l_i)\o \a(a)_{\lambda}(r_{j\mu}r_i)\o \a^{2}(l_j)|_{\mu=1\o 1\o \partial}|_{\lambda=-\partial^{\otimes^3}}\\
&=& -\sum_{i,j} \a(l_i)\o (a_\lambda r_j)_{\mu}r_i\o \a^{2}(l_j)|_{\lambda=-\partial^{\otimes^3},\mu=\partial^{\otimes^2}\o 1}\\
(D13)&=& \sum_{i,j} \a([l_{i\psi}r_j])\o \a(a_{\lambda}r_i)\o \a^{2}(l_j)|_{\mu=-1\o 1\o \partial-\partial\o 1\o 1}\\
&&+\a(r_{j\mu}l_i)\o \a(a)_{\lambda}\a(r_i)\o \a^{2}(l_j)|_{\mu=1\o 1\o \partial}|_{\lambda=-\partial^{\otimes^3}}\\
&=&\sum_{i,j} \a(l_{i\mu}r_j)\o \a(a_\lambda r_i) \o \a^2(l_j)|_{\lambda=-\partial^{\otimes^3}, \mu=-1\o 1\o \partial-\partial\o 1\o 1}\\
(D14)&=& \sum_{i,j}-\a(a_\lambda r_i)\o \a([l_{i\mu}r_j])\o \a^{2}(l_j)|_{\mu=-1\o \partial^{\otimes^2}}\\
&&-\a(a)_\lambda \a(r_i)\o \a(r_{j\mu}l_i)\o \a^{2}(l_j)|_{\mu=1\o 1\o \partial}|_{\lambda=-\partial^{\otimes^3}}\\
&=& -\sum_{i,j}\a(a_\lambda r_i)\o \a(l_{i\mu}r_j)\o \a^{2}(l_j)|_{\lambda=-\partial^{\otimes^3},\mu=-1\o \partial^{\otimes^2}}
\end{eqnarray*}
\begin{eqnarray*}
(D15)&=& \sum_{i,j}-\a^{2}(r_i)\o \a(a)_{\lambda}[l_{i\mu}r_j]\o \a^{2}(l_j)|_{\mu=\partial\o 1\o 1}+\a^2(r_i)\o (a_\lambda r_j)_{\mu}\a(l_i)\o \a^{2}(l_j)|_{\mu=\partial^{\otimes^2}\o 1}\\
&& -\a^2(r_i)\o \a(a)_\lambda (r_{j\mu}r_i)\o \a^{2}(l_j)|_{\mu=1\o1\o \partial}|_{\lambda=-\partial^{\otimes^3}}\\
&=& \a(r_i)\o [(a_\lambda r_j)_{\mu}r_i] \o  \a^2(l_j))|_{ \mu=- \partial^{\otimes^2}\o 1} -\a^2(r_i)\o [a_{\lambda}l_i]_{\mu}\a(r_j)\o \a^2(l_j) )|_{ \mu=-1\o \partial^{\otimes^2}}|_{\lambda=-\partial^{\otimes^3}}\\
(D16)&=& \sum_{i,j}  \a(a)_{\lambda}[l_{j\psi}r_i]\o \a^{2}(r_j)\o \a^{2}(l_i)|_{\mu=1\o \partial\o 1}   -(a_{\lambda}r_{j})_{\mu}l_i\o \a(r_i)\o \a^{2}(l_j)|_{\mu=\partial^{\otimes^2}\o 1}\\
&&+\a(a)_{\lambda}(r_{j\mu}l_i)\o \a^{2}(r_i)\o \a^{2}(l_j)|_{\mu=1\o 1\o \partial}|_{\lambda=-\partial^{\otimes^3}}\\
&=& \sum_{i,j}-[(a_\lambda r_j)_{\mu}l_i]\o \a(r_i)\o   \a^2(l_j))|_{ \mu=- \partial^{\otimes^2}\o 1}\\
&& +[a_{\lambda}l_i]_{\mu}\a(r_j)\o \a^2(r_i)\o\a^2(l_j)) |_{ \mu= -\partial\o 1\o 1-1\o 1\o \partial}|_{\lambda=-\partial^{\otimes^3}}\\
(D17)&=& \sum_{i,j} \a(r_{j\mu}r_i)\o \a(a_{\lambda} l_i)\o \a^{2}(l_j) |_{\mu=1\o 1\o \partial}-\a(a_{\lambda}l_i)\o \a(r_{j\mu}r_i)\o \a^{2}(l_j)  |_{\mu=1\o 1\o \partial}\\
&& +\a(a_{\lambda}l_i)\o \a(r_{j\mu}r_i)\o \a^{2}(l_j)  |_{\mu=1\o 1\o \partial}-\a(r_{j\mu}r_i)\o \a(a_{\lambda} l_i)\o \a^{2}(l_j) |_{\mu=1\o 1\o \partial}|_{\lambda=-\partial^{\otimes^3}}=0.
\end{eqnarray*}
Therefore, $(C1)=(D1)$. Similar to check $(C2)$ and $(C3)$. The proof is finished.   \hfill $\square$
\begin{proposition}
Let $(A, \a)$ be a  Hom-left-symmetric conformal  algebra which is free and finite as a $\mathbb{C}[\partial]$-module. The $\lambda$-product of $(A, \a)$ is obtained from the coalgebra $(A^{\ast c}, \a^{\ast}, \psi)$. Suppose that there is a  Hom-left-symmetric conformal  coalgebra structure $(A, \a, \varphi)$ where $\varphi$ is defined by (5.1). Then $\psi$ is a 1-cocycle of $\mathfrak{g}(A^{\ast c})$  associated to $(A^{\ast c}\o A^{\ast c}, \phi_{A^{\ast c}})$ if and only if $r$
\end{proposition}
\begin{center}
 {\bf ACKNOWLEDGEMENT}
 \end{center}

  The work of X. H. Zhang is  supported by  the Project Funded by China Postdoctoral Science Foundation (No. 2018M630768) and
the NSF of Shandong Province (No. ZR2016AQ03).
  The work of S. J. Guo is  supported by  the NSF of China (No. 11761017).
   The work of S. X. Wang is  supported by  the outstanding top-notch talent cultivation project of Anhui Province (No. gxfx2017123)
 and the Anhui Provincial Natural Science Foundation (1808085MA14).

\renewcommand{\refname}{REFERENCES}

\end{document}